\documentclass{article}
\usepackage{latexsym}
\font\Bbb=msbm10

\def\R{\hbox{\Bbb R}}
\def\eqref#1{(\ref{#1})}

\def\G{\Gamma}
\newtheorem{Th}{Theorem}
\newtheorem{Lemma}{Lemma}

\begin{document}

\title{Reconstruction of the connection or metric\\ from some partial information}

\author{J. Mike\v s, A. Van\v zurov\'a}

\date{\small
Fac.~Sci.,
Dept.~Algebra and Geometry,
Palack\'y University,\\
T\v r.~17.~listopoadu 12,
770 00 Olomouc,
Czech Republic,\\
mikes@inf.upol.cz, vanzurov@inf.upol.cz
}

\def\thefootnote{}
\footnotetext{\hspace*{-12pt} Supported by grant MSM 6198959214
of the Ministery of Education.}

\footnotetext{\hspace*{-12pt} \emph{Keywords: Riemannian manifold, linear connection,
metric, Fermi coordinates semigeodesic coordinates
}   }
\footnotetext{\hspace*{-12pt}
}

\def\thefootnote{\arabic{footnote}}
\def\proof{\noindent\textit{Proof. }}
\def\eproof{\hfill$\Box$}

\maketitle

\begin{abstract}

In \cite{G-A1}, in a neighborhood of a (positive definite) Riemannian space in which special, semigeodesic, coordinates are given, the metric tensor is calculated from its values on a suitable hypersurface and some of components of the curvature tensor of type $(1,3)$ in the coordinate domain.
Semigeodesic coordinates are a generalization of the well-known Fermi coordinates, that play an important role in mechanics and physics, are widely used in Minkowskian space (e.g. \cite{Ma},
\cite{Ma1}, \cite{A-Ma}), and in differential geometry of Riemannian spaces in  general.

In the present paper, we  consider a more general situation.
We introduce special pre-semigeodesic charts characterized both geometrically and in terms of the connection,
formulate a version of the 
Peano's-Picard's-Cauchy-like Theorem on existence and uniqueness of solutions
of the initial values problems for systems
of first-order ordinary differential equations. Then we use the apparatus in
a fixed pre-semigeodesic chart 
of a manifold 
equipped with the linear symmetric connection. Our aim is to
reconstruct, or construct, the symmetric linear connection in some neighborhood from the knowledge of the ``initial conditions":
the restriction of the connection to a fixed $(n-1)$-dimensional surface $S$ and some of the components of the curvature tensor $R$ in the ``volume"  (coordinate domain).
By 
analogous methods, we retrieve (or construct)
the metric tensor of type $(0,4)$ of a pseudo-Riemannian manifold in a domain
of semigeodesic coordinates from the known restriction of the metric to some
non-isotropic hypersurface and some of the components of the curvature tensor in the volume.
In comparison to the authors of \cite{G-A1}, we give shorter proofs of constructive character based on classical results
on first order ODEs (ordinary
differential equations).
\end{abstract}

\section{Introduction}

The problem of finding a Riemannian metric from this or that information is of interes from both theoretical and practical point of view. Papers by many authors are devoted to the possibility of finding the metric from the curvature tensor, \cite[p.~135-136]{Kow}, or prove existence of metrics with the prescribed Ricci tensor, \cite{deT},
\cite{Van} etc. In general, to solve the problem means to solve a relatively complicated non-linear system of partial differential equations, the coefficients of which are expressed through the components of the Riemannian curvature tensor. One possibility how to simplify the situation is to find a convenient coordinate system with respect to which the system of equations
is simplified considerably. Our aim is to present and use such preferable coordinates.

Recall the so-called Fermi coordinates, named after the Italian physicist Enrico Fermi,
\cite{F}, widely used in Riemannian geometry as well as in theoretical physics, \cite{A-Ma}, \cite{Ma}, \cite{Ma1} etc. Suppose $\gamma\colon I\to M$ is a geodesic on an $n$-dimensional Riemannian
manifold $M$, and $p$ a point on $\gamma$.
Then there exist local coor\-di\-na\-tes $(t,x^2,\dots ,x^n)$ around $p$ such that
for small $t$, $\gamma(t,0,\dots ,0)$ represents the geodesic near $p$;
on $\gamma$, the metric tensor is the Euclidean metric; and again (only) on $\gamma$,
all Christoffel symbols vanish (all the above properties are only valid on the distinguished geodesic). We will consider here a more general
si\-tu\-ation, namely pre-semigeodesic and semigeodesic coordinates, which bring, at the same time, special parametrization for all canonical geodesics in some tubular neighborhood.  The celebrated Fermi coordinates can be considered as a particular case.

\section{Pre-semigeodesic chart}

Let $(M,\nabla)$ be a (differentiable or smooth) $n$-dimensional manifold $M$ equipped with a symmetric linear connection $\nabla$.
Let ${\G}^h_{ij}$ denote components of the connection $\nabla$ in a fixed chart
$(U,\varphi=(x^1,\dots ,x^n))$ in $M$; $U\subseteq M$ open.

If in the chart $(U,(x^i))$ of $M$, ${\G}^h_{11}(x)=0$ is valid 
for all $h=1,\dots ,n$,
we say that  $(U,(x^i))$ is a \emph{pre-semigeodesic chart}\footnote{
Similar coordinates were used e.g.~in \cite{Pe}, and called there, in English translation, ``almost semigeodesic".}
related to the coordinate $x^1$
with respect to the connection $\nabla$ or,
that $x^1$ is a \emph{geodesic coordinate} in $U$. Obviously, it is quite natural to prefer the first coordinate,
and it means no loss of generality.

Let us give a geometric interpretation of the pre-semigeodesic charts.
Recall that the equations $\nabla_{\dot c}\dot c=0$ for canonically paramerized geodesics $c\colon I\to U$ of the connection
$\nabla\vert U$ in local coordinates read ($k=2,\dots ,n$)
\begin{equation}
\begin{array}{rl}
&\frac{d^2c^1}{ds^2}+\G^1_{11}\left(\frac{dc^1}{ds}\right)^2
+\sum_{j=2}^{n}\G^1_{1j}\frac{dc^1}{ds}\frac{dc^j}{ds}
+\sum_{i,j=2}^{n}\G^1_{ij}\frac{dc^i}{ds}\frac{dc^j}{ds}=0,\\
&\frac{d^2c^k}{ds^2}+\G^k_{11}\left(\frac{dc^1}{ds}\right)^2
+\sum_{j=2}^{n}\G^k_{1j}\frac{dc^1}{ds}\frac{dc^j}{ds}
+\sum_{i,j=2}^{n}\G^k_{ij}\frac{dc^i}{ds}\frac{dc^j}{ds}=0.
\end{array}
\label{1}
\end{equation}

\begin{Lemma}\label{L1}
The conditions $\G^h_{11}=0$, $h=1,\dots ,n$ are satisfied in $U$ if and only if the parametrized curves
\begin{equation}
c\colon I\to U,\qquad c(s)=(s,a_2,\dots ,a_n),\quad s\in I,\ a_i\in\R,\ \ i=2,\dots ,n
\label{2}
\end{equation}
are canonically parametrized geodesics of $\nabla\vert U$
($I$ is some interval, $a_k$ are suitable constants chosen so that $c(I)\subset U$).
\end{Lemma}

\proof
Let $\G^h_{11}=0$ hold  for $h=1,\dots ,n$. Then the local curves with parametrizations \eqref{2} satisfy
\begin{equation}
\frac{dc(s)}{ds}=\left(\frac{\partial }{\partial x^1}\right)_{c(s)}, \qquad \frac{d^2c(s)}{ds^2}=0,
\label{3}
\end{equation}
therefore are solutions to the system \eqref{1}.
Conversely, if the curves \eqref{2} are among solutions to \eqref{1}
then due to \eqref{3}, we get $\G^h_{11}=0$ from \eqref{1}.
\eproof

\smallskip
Hence the  pre-semigeodesic chart is fully characterized by the condition that the curves $x^1=s$, $x^i=\mbox{\rm const}$,
$i=2,\dots n$ belong to the geodesics of the given connection in the coordinate neighborhood.
The definition domain $U$ of such a chart is ``tubular", a tube along geodesics.

\section{Reconstruction of the connection}

Our aim is to show that a symmetric linear connection in a pre-semigeodesic coordinate domain $U$
(related to $x^1$)
can be uniquely constructed, or retrieved, in some subdomain of $U$, if we know the restriction
$\tilde\nabla=\nabla\vert S$ of the connection to the surface $S$
defined by $x^1=0$ and the prescribed components $R^h_{i1k}$
of the curvature tensor in the given tubular domain $U$.

First let us modify, for our purpose, the Theorem on existence and uniqueness of solutions of systems
of ODEs.

%
In $\R^n$ with standard coordinates $(x^1,x^2,\dots,x^n)$, let us identify the li\-near subspace (hypersurface) characterized
by $x^1=0$ with $\R^{n-1}$, i.e.~$(\tilde x)=(x^2,\dots,x^n)$ are standard coordinates in ${\R}^{n-1}$. Let
${\mathcal J}=(0,1)$ be the open unit interval
and denote by $K_{m}={\mathcal J}^{m}$ the open standard $m$-cube.
Denote
$$D_{n}(\delta)=\{x=(x^1,\dots, x^n)\in\R^{n}\mid 0\le x^1\le\delta,\ 0<x^i<1, \ i=2,\dots ,n\}
.$$
The open $(n-1)$-cube $K_{n-1}={\mathcal J}^{n-1}$, viewed as
$$K_{n-1}=\{\tilde x=(x^2,\dots, x^n)\in\R^{n-1}\mid 0<x^i<1, \
i=2,\dots ,n\}\subset\R^{n-1},$$
can be identified with a hypersurface $S$ in $D_{n}(\delta)$ determined by $x^1=0$.

\begin{Th}\label{T1}
Let $S$ be a hypersurface in $D_{n}(\delta)$ defined by $x^1=0$.
Let $\tilde\nabla$ be a symmetric linear connection in $S$
(of the class at least $C^2$)
with the components ${\tilde\G}^h_{ij}$ and the curvature tensor $\tilde R$,
and let $A^h_{ij}$ be functions in $D_{n}(\delta)$
(at least $C^0$ in $D_{n}(\delta)$)
and such that each
$A^h_{1k}$ is at least $C^1$ in each of the variables $x^2,\dots ,x^n$ and
at least $C^0$ in $x^1$. Moreover, let the conditions
${\tilde R}^h_{i1k}=A^h_{ik}$ hold in $S$.
Then there is a real number $\hat\delta$, $0<\hat\delta\le\delta$ and
a unique symmetric linear connection $\nabla$ in some $D_{n}(\hat\delta)$
with components satisfying $\G^h_{11}=0$, $h=1,\dots ,n$, such that
$\nabla\vert S=\tilde\nabla$ and
${R}^h_{j1k}={A}^h_{jk}$ in $D_{n}(\hat\delta)$ for $j=1,\dots ,n$.
\end{Th}

\proof
The components of the curvature tensor are related to the components of the connection by the classical formula
\begin{equation}
R^h_{ijk}={\partial}_j{\G}^h_{ik}-{\partial}_k{\G}^h_{ij}+{\G}^m_{ik}{\G}^h_{mj}-{\G}^m_{ij}{\G}^h_{mk}.
\label{4}
\end{equation}
Particularly under the assumption ${\G}^h_{11}$, setting $i=j=1$ we get
\begin{equation}
\frac{\partial}{\partial x^1}{\G}^h_{1k}+\sum_{m}{\G}^m_{1k}{\G}^h_{1m}-R^h_{11k}=0.
\label{5}
\end{equation}
Let us write the system as
\begin{equation}
\frac{\partial}{\partial x^1}{\G}^h_{1k}=-\sum_{m}{\G}^m_{1k}{\G}^h_{1m}+A^h_{1k}.
\label{5a}
\end{equation}
We can view (\ref{5a}) as a system of ODEs 
of one variable $x^1$,
while the remanining coordinates $(\tilde x)=(x^2,\dots ,x^n)\in K_{n-1}=S$
are considered as parameters.
Given the initial data
${\tilde\G}^h_{1k}(\tilde x)$, $\tilde x\in S$,
there exists $\delta_1$, $0<\delta_1\le \delta$ and there are uniquely determined
functions ${\G}^h_{1k}(x^1,\dots,x^n)$
of the class at least $C^1$
on the domain $D_{n}(\delta_1)$ such that
\begin{equation}
{\G}^h_{1k}(0,\tilde x)={\tilde\G}^h_{1k}(\tilde x),\qquad
\tilde x\in S.
\label{6}
\end{equation}
Now setting $j=1$ in \eqref{4} we get for the indices $i=2,\dots ,n$ the system
\begin{equation}
\frac{\partial}{\partial x^1}{\G}^h_{ik}+\sum_{m}{\G}^m_{ik}{\G}^h_{m1}-
\frac{\partial}{\partial x^k}{\G}^h_{i1}
-R^h_{i1k}=0
\label{7}
\end{equation}
which we rewrite as
\begin{equation}
\frac{\partial}{\partial x^1}{\G}^h_{ik}=-\sum_{m}{\G}^m_{ik}{\G}^h_{m1}+
\frac{\partial}{\partial x^k}{\G}^h_{i1}
+A^h_{ik}=0.
\label{7a}
\end{equation}
Substituting the obtained functions ${\tilde\G}^h_{1k}$ we find
that according to the existence and uniqueness theorem on systems of ODEs
there is $\hat\delta$, $0<\hat\delta\le \delta_1$ and
there are uniquely determined
functions ${\G}^h_{ik}(x^1,\dots,x^n)$
of the class at least $C^1$
on the domain $D_{n}(\hat\delta)$ such that the initial conditions
\begin{equation}
{\G}^h_{ik}(0,\tilde x)={\tilde\G}^h_{ik}(\tilde x),\qquad
\tilde x\in S
\label{8}
\end{equation}
are satisfied. Moreover, we can easily see that due to \eqref{5}, \eqref{5a}, \eqref{7}, \eqref{7a},
\begin{equation}
{R}^h_{i1k}( x)={A}^h_{ik}( x),\qquad
x\in  D_{n}(\hat\delta),\ \ i=1,\dots ,n
\label{9}
\end{equation}
holds as required.
\eproof

\smallskip
As a consequence, if we use eventual
prolongation of the solution, we obtain:

\begin{Th}\label{T2}
Let $(U,\varphi=(x^1,\dots ,x^n))$ be a chart in $M$.
Let $S\subset U$ be a submanifold in $U$ defined by $x^1=0$.
Let $\tilde\nabla$ be a symmetric linear connection in $S$ of the class at least $C^2$
with the components ${\tilde\G}^h_{ij}$ and the curvature tensor $\tilde R$,
and let $A^h_{ij}$ be functions in $U$ such that ${\tilde R}^h_{i1k}=A^h_{ik}$ in $S$,
$A^h_{ik}$, $i=2,\dots ,n$ are of the class at least $C^0$,  $A^h_{1k}$ are continuous in
$x^1$, and $A^h_{1k}$ are at least $C^1$ in the remaining variables
$x^2,\dots ,x^n$.
Then there is a unique symmetric linear connection $\nabla$ in $U$ with components satisfying $\G^h_{11}=0$
for $h=1,\dots ,n$ (i.e.~the given chart is pre-semigeodesic
w.r.t.~$\nabla$) such that $\nabla\vert S=\tilde\nabla$, and ${R}^h_{j1k}={A}^h_{jk}$ in $U$.
\end{Th}

\section{Reconstruction of a metric}

\subsection{Semigeodesic coordinates}

For our purpose, we say that a chart $(U,(x^i))$ of a pseudo-Riemannian ma\-ni\-fold $(M,g)$ is \emph{semigeodesic}
(or that $(x^i)$ are semigeodesic coordinates) if in this chart,
the metric tensor has  the coordinate expression
\begin{equation}
g=edx^1\otimes dx^1+{\tilde g}_{ij}(x^1,
x^2,\dots ,x^n) dx^i\otimes dx^j,\quad i,j=2,\dots ,n
\label{m}
\end{equation}
where $e=\pm 1$ (the plus or minus sign is connected with the square of the integral of the tangent vector to the $x^1$-coordinate line).

The geometric interpretation is as follows, \cite[p.~55]{Pe}.

\begin{Lemma}\label{L2}
Local coordinates $(x^i)$ in a pseudo-Riemannian ma\-ni\-fold are semigeodesic if and only if the 1-net of $x^1$-coordinate lines is formed by
arclength parametrized geodesics which are orthogonal to a non-isotropic hypersurface defined by $x^1=\mbox{\rm const}$.
\end{Lemma}

Note that coordinate hyperplanes defined by $x^j=\mbox{\rm const}$ are orthogonal to the distinguished system of geodesics.
Obviously, semigeodesic coordinates are pre-semigeodesic.

Semigeodesic coordinates can be introduced in a sufficiently small neighborhood of any point of an arbitrary (positive)
Riemannian manifold, and is fully characterized by the coordinate form of the metric:
\begin{equation}
g_{ij}=dx^1\otimes dx^1+{\tilde g}_{ij}(x^1,
x^2,\dots ,x^n) dx^i\otimes dx^j,\quad i,j=2,\dots ,n.
\label{10}
\end{equation}
E.g.~on a cylinder, semigeodesic coordinates can be introduced globally.

Advantages of such coordinates are known since Gauss (\cite[p.~201]{KR}, ``Geod\" a\-ti\-sche Parallel\-koordinaten"), and are widely used in the two-dimensional case, particularly in applications, \cite{MKY} and the references therein, \cite{Ve} etc.
Note that geodesic polar coordinates (``Geod\" atische Polarkoordinaten,"
\cite[pp.~197-204]{KR}) can be  interpreted as a ``limit case" of semigeodesic coordinates (all geodesic coordinate lines $\phi=x^2=\mbox{\rm const}$ pass through one point called the pole, cor\-res\-ponding to
$r=x^1=0$, while $r=x^1=\mbox{\rm const}$ are the geodesic circles).


\subsection{Reconstruction of a metric in the semigeodesic coordinates}

\begin{Th}\label{T3}
Let $a_{ij}$ be (at least) continuous functions in $D_{n}(\delta)$, let ${\tilde g}_{ij}$ be functions of the class
(at least) $C^2$ in $K_{n-1}$ and
${\tilde G}_{ij}$ functions of the class (at least) $C^1$ in $K_{n-1}$, $i,j=2,\dots ,n$ such that the matrices
$({\tilde g}_{ij})$ and $({\tilde G}_{ij})$ are symmetric\footnote{${\tilde g}_{ji}={\tilde g}_{ij}$,\
${\tilde G}_{ji}={\tilde G}_{ij}$
}
and $\det({\tilde g}_{ij})\ne 0$ in $K_{n-1}$.
Fix an element $e\in\{-1,1\}$.
Then there is $\hat\delta$, $0<\hat\delta\le \delta$
and there exists exactly one non-degenerate metric tensor\footnote{
$\det({g}_{ij})\ne 0$ on $D_{n}(\hat\delta)$}
$g$ of the class (at least) $C^2$ in $D_{n}(\hat\delta)$ with components
$g_{11}=e$, $g_{1j}=0$, $j=2,\dots ,n$
such that for $i,j=2,\dots ,n$,
\begin{equation}
g_{ij}(0,\tilde x)={\tilde g}_{ij}(\tilde x), \quad \frac{\partial^+}{\partial x^1}\,g_{ij}(0,\tilde x)={\tilde G}_{ij}(\tilde x),
\quad\tilde x\in K_{n-1} 
\label{10a}
\end{equation}
where $\frac{\partial^+}{\partial x^1}$ means the partial derivative from the right, and
\begin{equation}
a_{ij}( x)={R}_{1ij1}(x), \qquad x\in D_{n}(\hat\delta).
\label{11}
\end{equation}
\end{Th}

\proof
The components of the curvature tensor $R$ (in type $(0,4)$) of the semi-Riemannian manifold $V_n=(M,g)$ are related to the components of the metric by
\begin{equation}
R_{hijk}=\frac{1}{2}({\partial}_{ij}g_{hk}
+{\partial}_{hk}g_{ij}
-{\partial}_{ik}g_{hj}
-{\partial}_{ij}g_{hk}
)
+g^{rs}({\G}_{hkr}{\G}_{ijs}-{\G}_{hjr}{\G}_{kjs})
\label{11a}
\end{equation}
where ${\G}_{ijk}=\frac{1}{2}({\partial}_{i}g_{jk}
+{\partial}_{j}g_{ik}
-{\partial}_{k}g_{ij})$ are Christoffel symbols of the first type in $V_n$, and $g^{rs}$ are components of the dual tensor to $g$.
Hence $g^{ij}$ are functions rational in components $g_{ij}$ of the metric\footnote{
$g^{ij}=1\slash \det(g_{ij})\cdot  A_{ji}$ where $A_{ji}$
is the algebraic complement of the matrix element $g_{ji}$}.

Setting $h=k=1$ and using the assumptions $g_{11}=e$, $g_{1j}=0$
we obtain from \eqref{11a}
\begin{equation}
R_{1ij1}=\frac{1}{2}{\partial}_{11}g_{ij}
-\frac{1}{4}{g}^{rs}{\partial}_{1}g_{ir}{\partial}_{1}g_{js}.
\label{12}
\end{equation}
Here we can suppose that the indices satisfy $i,j,r,s>1$.
Let us substitute
\begin{equation}
G_{ij}={\partial}_{1}g_{ij}.
\label{13}
\end{equation}
Then we can write \eqref{12} as
\begin{equation}
R_{1ij1}=\frac{1}{2}{\partial}_{1}G_{ij}
-\frac{1}{4}{g}^{rs}G_{ir}G_{js}.
\label{14}
\end{equation}
Denote $R_{1ij1}(x)=a_{ij}(x)$.
Hence we should solve the system
\begin{equation}
\displaystyle
\begin{array}{rl}
{\partial}_{1}g_{ij}&=G_{ij},\\
&\\
{\partial}_{1}G_{ij}&=\frac{1}{2}{g}^{rs}G_{ir}G_{js}+2a_{ij}
\end{array}
\label{15}
\end{equation}
under the initial values
\begin{equation}
g_{ij}(0,\tilde x)={\tilde g}_{ij}(\tilde x), \quad \frac{\partial^+}{\partial x^1}\,g_{ij}(0,\tilde x)={\tilde G}_{ij}(\tilde x),\quad \tilde x\in K_{n-1},
\ i,j=2,\dots ,n.
\label{16}
\end{equation}
Note that since the determinant as well as the algebraic complements are continuous functions in the entries $g_{ij}$, and we demand
$\det({\tilde g}_{ij})(0,\tilde x) =\det({\tilde g}_{ij})(\tilde x)\ne 0$, it is guaranteed that ${g}^{rs}$ will be well-defined and well-behaved functions of ${g}^{ij}$, similarly as in \cite{G-A1}.
So \eqref{15} can be considered as a system of first-order ordinary differential equations
in the variable $x^1$ for the unknown functions $g_{ij}$ and $G_{ij}$ with the initial values
\eqref{16}; the remaining coordinates $x^2,\dots ,x^n\in K_{n-1}$ are supposed to be parameters.
The right sides in (\ref{15}) satisfy the conditions of the existence and uniqueness theorem \cite[p.~263]{D-G} in the domain
$D_n(\tilde\delta)$ and have continuous derivatives with respect to $g_{ij}$ and $G_{ij}$. The initial value problem
(\ref{15}) and (\ref{16}) has precisely one solution $g_{ij}(x)$. The functions $g_{ij}$ are components of a metric tensor in
$D_n(\tilde\delta)$, and comparing (\ref{15}) and (\ref{14}) we find easily that the components of its curvature tensor satisfy
$R_{1ij1}(x)=a_{ij}(x)$ as required.
\eproof

\medskip

Since the matrices $(g_{ij})$ and $(G_{ij})$ are symmetric we may assume $i\le j$ 
in (\ref{15}) and (\ref{16}).

As a consequence, we get

\begin{Th}\label{T4}
Let $a_{ij}$ be continuous functions in some coordinate neighborhood $U$,
${\tilde g}_{ij}$ $C^2$-functions in $\tilde S=U\cap S$ where $S$ is the hypersurface $S\colon x^1=0$ in $R^n$, and
${\tilde G}_{ij}$ $C^1$-functions in $\tilde S$, $i,j=2,\dots ,n$
such that the matrices
$({\tilde g}_{ij})$ and $({\tilde G}_{ij})$ are symmetric
and $\det({\tilde g}_{ij})\ne 0$ in $\tilde S$.
Fix an element $e\in\{-1,1\}$.
Then there is $\hat\delta>0$
and there exists precisely one
non-degenerate metric tensor $g$,
$\det({g}_{ij})\ne 0$, of the class $C^2$ in $\tilde U=\langle -\hat\delta, \hat\delta\rangle \times\tilde S $
with components
$g_{11}=e$, $g_{1j}=0$, $j=2,\dots ,n$ (i.e.~$\tilde U$ is semigeodesic)
such that for $i,j=2,\dots ,n$,
\begin{equation}
g_{ij}(0,\tilde x)={\tilde g}_{ij}(\tilde x), \quad {\frac{\partial^{+}}{\partial x^1}}\,g_{ij}(0,\tilde x)={\tilde G}_{ij}(\tilde x),\quad\tilde x\in \tilde S
\label{10b}
\end{equation}
and
\begin{equation}
a_{ij}( x)={R}_{1ij1}(x), \qquad x\in \tilde U.
\label{11b}
\end{equation}
\end{Th}

\medskip

Provided $a_{ij}(x)=R_{1ij1}(x)$ the solution of the system (\ref{15}) answers the problem of finding the metrics with the prescribed
components $R_{1ij1}(x)$ of the $(0,4)$-Riemannian curvature tensor. Substituting the obtained components of metric we get the relationship
to the components of the $(1,3)$-curvature as follows:

\begin{equation}
R_{1ij1}=e R^1_{ij1}=-eR^1_{i1j}=g_{im}R^m_{11j}=-g_{im}R^m_{1j1}.
\label{17}
\end{equation}
Hence our results generalize the results of \cite{G-A} and \cite{G-A1}.

\end{document}